\documentclass[12pt]{amsart}
\usepackage{amscd,amssymb,graphics}
\newcommand\mathscr{\mathcal}

\newtheorem{theorem}{Theorem}
\newtheorem{corollary}[theorem]{Corollary} 
\newtheorem{lemma}[theorem]{Lemma}

\theoremstyle{definition}
\newtheorem{definition}[theorem]{Definition}

\theoremstyle{remark}

\renewcommand{\Bbb}{\mathbb}

\def\norm#1{\left\Vert#1\right\Vert}

\def\I {{\Bbb I}}
\def\C {{\Bbb C}}
\def\N{{\Bbb N}}

\def\Id{{\mathbb{Id}}}
\def\U{{\Bbb U}}

\def\Aut{{\mathrm{Aut}}\,}

\def\diam{{\mathrm{diam}}}

\def\H{{\mathcal H}}

\def\QED{\nobreak\quad\ifmmode\roman{Q.E.D.}\else{\rm Q.E.D.}\fi}

\def\ss{\mathbb{S}}

\def\G{\Gamma}

\def\sbs{\subset}

\def\e{\varepsilon}

\def\ti{\times}
\def\obr{^{-1}}

\def\Rep{\mathrm{Rep}\,}
\def\Iso{\mathrm{Iso\,}}

\begin{document}

\title[Representations by isometries]
{Representations of residually finite groups by 
isometries of the Urysohn space}

\author[V.G. Pestov]{Vladimir G. Pestov}

\address{Department of Mathematics and Statistics, 
University of Ottawa, 585 King Edward Ave., Ottawa, Ontario, Canada K1N 6N5}

\email{vpest283@uottawa.ca}

\author[V.V. Uspenskij]{Vladimir V. Uspenskij}

\address{Department of Mathematics, 321 Morton Hall, Ohio
University, Athens, Ohio 45701, USA}

\email{uspensk@math.ohiou.edu}

\thanks{Research by the first named author supported by NSERC discovery grant 
261450-2003 and University of Ottawa internal grants.}



\begin{abstract} 
As a consequence of Kirchberg's work,
Connes' Embedding Conjecture is equivalent to the property that 
every homomorphism of the group $F_\infty\times F_\infty$ into the unitary group
$U(\ell^2)$ with the strong topology is pointwise approximated by homomorphisms 
with a precompact range. In this form, the property (which we
call Kirchberg's property) makes sense for an arbitrary topological group. We establish the validity of the Kirchberg property for the isometry
group $\Iso(\U)$ of the universal Urysohn metric space $\U$ as a consequence
of a stronger result: every representation of a residually finite group 
by isometries of $\U$ can be pointwise
approximated by representations with a finite range. This brings up the natural
question of which other concrete infinite-dimensional groups satisfy the
Kirchberg property.  
\end{abstract}

\maketitle

\section{Introduction} 
Our motivation comes from theory of operator algebras. Recall that a 
$C^\ast$-algebra is {\em residually finite dimensional} if
it admits a separating family of $\ast$-homomorphisms into finite-dimensional
matrix algebras. 
The question of describing those countable discrete
groups $\Gamma$ for which
the full group $C^\ast$-algebra, $C^\ast(\Gamma)$, is residually finite-dimensional is of considerable interest.
For instance, Kirchberg's work \cite{kirchberg} implies that the validity of Connes Embedding Conjecture is equivalent to the
residual finite dimensionality of the algebra $C^\ast(F_2\times F_2)$.

A group $\G$ is {\em residually finite} if it admits a separating family
of homomorphisms to finite groups, or, equivalently, if the intersection
of all subgroups of finite index is reduced to the unity.
In view of Malcev's theorem stating that finitely generated
subgroups of linear groups are residually finite \cite{malcev}, 
for a finitely generated group $\Gamma$ residual finite-dimensionality of
$C^\ast(\Gamma)$ implies residual finiteness of $\Gamma$. However,
as shown by Bekka \cite{bekka}, this necessary condition on $\Gamma$ is not sufficient. 

It is easily seen that the $C^\ast$-algebra $C^\ast(\Gamma)$ of a countable
group $\Gamma$ is residually finite dimensional if and only
if every unitary representation of $\Gamma$ in $\ell^2$
can be approximated by representations of $\Gamma$ which factor through
representations of compact groups. Here the topology on the set of
representations $\Rep(\Gamma,\ell^2)$ is that induced from the product topology 
on $U(\ell_2)^\Gamma$, where the unitary group $U(\ell^2)$ is equipped with the
strong operator topology. 
Since the group $U(\ell^2)$ 
can be viewed as a group of isometries
of the unit sphere $\ss^{\infty}$ in $\ell^2$, a way to
a better understanding of the problem could be to enlarge the settting and consider representations of
$\Gamma$ by isometries of other metric spaces.                          
Other than the sphere $\ss^{\infty}$, a very natural candidate to consider 
is the universal Urysohn metric space $\U$ \cite{U25, U27, U72}, 
the object of considerable attention in recent years 
\cite{Kat,Usp90,Ver98,Gr,GK,P05}. Just like the sphere, the Urysohn space
is an ultrahomogeneous complete separable metric space
(that is, every isometry between two finite subspaces extends to a global
self-isometry of the entire space), but unlike the sphere,
$\U$ contains an isometric copy of every separable metric space.

For a topological group $G$ and a discrete group $\Gamma$, denote by 
$\Rep (\Gamma, G)$ the space of all homomorphisms from $\Gamma$ to $G$, considered as a subspace of $G^{\Gamma}$ with the product topology.
A representation of a group $\Gamma$ by isometries of a metric space $X$ is
just a homomorphism $\Gamma\to\Iso(X)$, where we equip the latter group with
the topology of pointwise convergence on $X$ (it coincides with the 
compact-open topology). We will denote $\Rep(\Gamma,\Iso(X))$ simply by
$\Rep(\Gamma,X)$.

\begin{theorem} 
\label{th:app}    
Let $\Gamma$ be a residually finite group. Every representation of
$\Gamma$ by isometries of the Urysohn space $\U$ is approximated 
by representations with a finite range. In other words,
representations of $\Gamma$
with a finite range are everywhere dense in the space $\Rep(\Gamma,\U)$.
\end{theorem}             

Our method of proof further refines a construction from \cite{P02,P06}. 
Note again that for unitary representations an analogue of Theorem \ref{th:app} does not hold, as follows from above mentioned Bekka's result from \cite{bekka} together with Corollary \ref{c:el} below.

Recall that
the \textit{lower Vietoris topology} on the set ${\mathscr F}(X)$
of all closed subsets of a topological space $X$ is determined by the
basic sets
\[L(V_1,V_2,\ldots,V_n) = \{F\in {\mathscr F}(X)\colon \forall i=1,2,\ldots,n,~~
F\cap V_i\neq\emptyset\},\]
where $V_i$, $i=1,2,\ldots,n$ are open subsets of $X$. 

Consider the following three properties of a topological group $G$:

\bigskip

\noindent 
($\circledast_1$)
if $A$ and $B$ are finite subsets of $G$ and every element of $A$ commutes
with every element of $B$, then there exist finite subsets $A'$ and $B'$ of $G$
that are
arbitrarily close to $A$ and $B$, respectively, 
such that every element of $A'$
commutes with every element of $B'$ and the subgroups of $G$ generated by $A'$ and $B'$ are relatively compact (equivalently, the subgroup generated by $A'\cup B'$ 
is relatively compact). 

\bigskip                                                     

Here ``arbitrarily close'' means the following: if $A=\{x_1,\dots, x_n\}$ and
$Ox_1,\dots, Ox_n$ are neighbourhoods of $x_1,\dots, x_n$, then $A'=\{x_1',\dots, x_n'\}$,
where $x_i'\in Ox_i$, and similarly for $B'$. 
If $G$ is a metric group, then one can simply talk of closeness with regard
to the Hausdorff metric. ``Relatively compact'' means ``with a compact
closure.''

\bigskip

\noindent 
($\circledast_2$)
For every pair $H_1$, $H_2$ of closed topological subgroups
of $G$ where every element of $H_1$ commutes with every element of
$H_2$, there are nets $K_{1,\alpha}$,
$K_{2,\alpha}$ of compact subgroups of $G$ such that 
\begin{itemize}
\item Every element of $K_{1,\alpha}$ commutes with every element of
$K_{2,\alpha}$,
and
\item $K_{i,\alpha}$ converges to $H_i$ in the lower Vietoris topology,
$i=1,2$.
\end{itemize}

\bigskip

\noindent 
($\circledast_3$)
Let $F_\infty$ be the free group on countably many generators, 
$\G=F_\infty\times F_\infty$. Then the homomorphisms $h:\G\to G$ with a relatively
compact range form a dense subset of the space $\Rep(\G,G)$. 

\bigskip

The reader will easily verify that the properties 
$\circledast_1$, $\circledast_2$,
$\circledast_3$ are equivalent to each other. We confine ourselves by the following
remark: every homomorphism of $\G=F_\infty\times F_\infty$ to $G$ gives rise to
a pair of commuting subgroups of $G$, and every pair of commuting finitely-generated
subgroups of $G$ can be obtained in such a way.

\begin{definition} 
A topological group $G$ satisfies the
{\em Kirchberg property} if it satisfies the equivalent properties
$\circledast_1$, $\circledast_2$, $\circledast_3$ considered above.
\end{definition}

By force of Kirchberg's work, 
Connes' Embedding Conjecture is equivalent to 
the Kirchberg property for the unitary 
group $G=U(\ell^2)$ with the strong topology. In our view,
it makes sense to investigate the validity of the Kirchberg property for
other infinite-dimensional topological groups of importance.
In particular,
it follows from Theorem \ref{th:app} that the isometry group of the Urysohn
space, which forms an infinite-dimensional group of considerable current
interest, satisfies the property in question.

\begin{corollary}
\label{c:kirch}
The isometry group $\Iso(\U)$ of the Urysohn metric space $\U$,
equipped with the natural Polish topology (topology of simple convergence), satisfies the Kirchberg property.
\end{corollary}

Our paper is organized as follows. 
The proof of our main Theorem \ref{th:app}, as well as that of Corollary \ref{c:kirch}, are contained in Section
\ref{s:main}, while the necessary prerequisites about the Urysohn space and
its group of isometries are collected in Section \ref{s:prereq}.
Finally, in Section \ref{s:cec} we discuss in some detail Connes' Embedding
Conjecture and the proposed Kirchberg property of topological groups.

\section{\label{s:prereq}Prerequisites on the Urysohn space}

Recall that the universal Urysohn metric space $\U$ 
\cite{Usp90,Ver98,Gr,GK,P05} is characterized by the following properties: 

(1) $\U$ is a complete separable metric space;

(2) $\U$ contains an isometric copy of every separable metric space;

(3) every isometry between two finite metric subspaces
of $\U$ extends to a global isometry of $\U$ onto itself.

Alternatively, $\U$ is the only complete separable metric space that is
{\em finitely injective:} for every finite metric spaces $K\sbs L$ every
distance-preserving map $K\to\U$ has a distance-preserving extension $L\to\U$. 
As a topological space, $\U$ is the same as the Hilbert space $\ell^2$ 
\cite{Usp04}.
The group $P=\Iso(U)$ is a universal topological group with a countable base:
every topological group with a countable base is isomorphic to a topological
subgroup of $P$ \cite{Usp90}, \cite{Usp98}, \cite{Usp02}. 
In this section we consider a few lemmas that will be used in Section 3.

\begin{lemma}
\label{l:Pest25}
Let $G$ be a residually finite group, $d$ a left-invariant pseudometric
on $G$, $K$ a finite subset of $G$. Suppose that the restriction of $d$
to $K$ is a metric. Then there exist a subgroup
$H\sbs G$ of finite index and a $G$-invariant pseudometric on 
the (finite) set $G/H$ such that the restriction of the quotient map
$G\to G/H$ to $K$ is distance-preserving.
\end{lemma}

This is essentially Lemma 2.5 of \cite{P06}. For the reader's 
convenience we give a proof. 

\begin{proof}
A {\em seminorm} on a group is a non-negative function $p$ satisfying the following:
(1) $p(1)=0$; (2) $p(xy)\le p(x)+p(y)$; (3) $p(x\obr)=p(x)$. 
There is a one-to-one correspondence
between left-invariant pseudometics and seminorms: given a seminorm $p$, consider
the pseudometric $\mu$ defined by $\mu(x,y)=p(x\obr y)$; given $\mu$, we recover $p$ by
$p(x)=\mu(x,1)$. 

The values of pseudometrics and seminorms are (finite) non-negative real numbers,
but in this proof it will be convenient to allow the value $+\infty$ as well.
In that case we use the terms {\em pseudometric$_*$} and {\em seminorm$_*$.}

Let $p$ be the seminorm corresponding to $d$, $p(x)=d(x, 1)$. 
Let $L=K\obr K\sbs G$. Consider the greatest seminorm$_*$ $q$
on $G$ that agrees with $p$ on $L$. This seminorm$_*$ is defined by
\[
q(x)=\inf \sum_{i=1}^n p(y_i),
\]
where the infimum is taken over all representations of the form $x=y_1\dots y_n$,
$y_i\in L$, $1\le i\le n$. (The infimum of the empty set is $+\infty$, thus $q$
is infinite outside the subgroup generated by $K$.) As $p$ is strictly positive
on $L\setminus\{1\}$, for every $C>0$ the set $G_C=\{x\in G:q(x)\le C\}$ is finite.
Let $M=\max\{p(x):x\in L\}$ and $C=2M$.  
Since $G$ is residually finite, there exists a normal
subgroup $H\triangleleft\, G$ of finite index such that $H\cap G_C=\{1\}$.   
Consider the quotient seminorm$_*$ $\dot q$ on $G/H$, defined by 
\[
\dot q(xH)=\inf\{q(y): y\in xH\}.
\]
We claim that $\dot q(xH)=p(x)$ for every $x\in L$. Indeed, otherwise there
exist $h\in H$ and a representation $xh=y_1\dots y_n$ with $y_1,\dots, y_n\in L$ 
and $\sum p(y_i)< p(x)$. But then $h\ne 1$ and 
$q(h)=q(x\obr y_1\dots y_n)\le p(x)+\sum p(y_i)
<2p(x)\le C$, hence $h\in H\cap G_C$, in contradiction with our choice of $H$. 

Equip $G/H$ with the pseudometric$_*$ $\lambda$ corresponding to $\dot q$. 
Then the quotient map $G\to G/H$ is distance preserving on $K$: if $x,y\in K$,
then $x\obr y\in L$, so $\lambda(xH,yH)=\dot q(x\obr yH)=p(x\obr y)=d(x,y)$.
If $\lambda$ takes infinite values, replace it by $\inf(\lambda, M)$.
\end{proof}

\begin{lemma}
\label{l:appr}
Let $\mu$ and $\lambda$ be two pseudometrics on a finite set $K$. 
Let $R$ be the diameter of the space $(K,\lambda)$. For every
distance-preserving map $i: (K,\mu)\to \U$ there exists a distance preserving
map $j: (K, \mu+\lambda)\to \U$ such that $d(i(x),j(x))\le R$ for every $x\in K$.
Moreover, if $a\in K$ is given, we may require that $j(a)=i(a)$. 
\end{lemma}

\begin{proof}
Consider the product $K\ti K$ with the pseudometric $\nu$ defined by 
\[
\nu((x,y),(x',y'))=\mu(x,x')+\lambda(y,y')
\]
and the embeddings $i',j':K\to K\ti K$ defined by
$i'(x)=(x,a)$ and $j'(x)=(x,x)$.
Write $i$ as $i=g\circ i'$, where $g:K\ti\{a\}\to U$ is distance-preserving.
Since $U$ is finitely injective, there exists a distance-preserving map 
$h:K\ti K\to U$ extending $g$. Set $j=h\circ j'$. Then $j(a)=h(a,a)=i(a)$ and
$d(i(x),j(x))=d(h(i'(x)),h(j'(x)))=
\nu(i'(x),j'(x))=\lambda(x,a)\le R$ for every $x\in K$. 
\end{proof}

Recall that the group $\Iso(X)$ of isometries of a metric space is equipped
with the topology of pointwise convergence on $X$, or, which is equivalent, the
compact-open topology. If $X$ is complete separable, the group $\Iso(X)$ is
Polish.

Following \cite{Usp98}, let us say that a subspace $X\sbs \U$ is 
{\em g-embedded} if there exists a homomorphism $k:\Iso(X)\to\Iso(\U)$
of topological groups 
such that for every $h\in \Iso(X)$ the isometry $k(h)\in \Iso(\U)$ is
an extension of $h$. 

\begin{lemma}
\label{l:g-embed}
Every finite subset of $U$ is $g$-embedded.
\end{lemma}

\begin{proof}
Every separable metric space admits a $g$-embedding in $\U$ 
\cite{Usp90, Usp98, Usp02}. Since any two embeddings of a finite metric
space in $\U$ are conjugate by an isometry of $U$, the lemma follows.
\end{proof}

The same argument shows that ``finite'' can be replaced by ``compact'' in 
the lemma.

\section{\label{s:main}Proofs of the main results}

\subsection*{Proof of Theorem \ref{th:app}}
We denote the topological group $\Iso(\U)$ by $P$. Let $\G$ be a residually
finite group.
Let $f:\G\to P$ be a homomorphism, and let $O$ be a neighborhood of $f$
in $\Rep(\G,P)$. We must find $f'\in O$ with a finite range.

We may assume that $O$ has the following form: there are finite sets 
$K\sbs \G$ and $A\sbs \U$ and $\e>0$ 
such that any $h\in \Rep(\G,P)$ is in $O$ if and only
if 
$h(g)x$ is $\varepsilon$-close to $f(g)x$ for every $g\in K$ and $x\in A$.
(We say that two points are {\em $\e$-close} 
if the distance between them is $\le\e$.)
We may also assume that $K$ contains the neutral element $1\in\G$. 

{\em Case 1.} Suppose first that $A$ is a singleton: $A=\{p\}$ for some $p\in \U$.
Let $d$ be the metric on $\U$. Consider the left-invariant pseudometric $\mu$
on $\G$ defined by $\mu(g_1,g_2)=d(f(g_1)p, f(g_2)p)$. Let $\lambda$ be the discrete
metric on $\G$ such that $\lambda(g,h)=1$ for any pair of distinct elements $g,h\in \G$.
Let $\rho=\mu+\e\lambda$. This is a left-invariant metric on $\G$.
According to 
Lemma~\ref{l:Pest25}, there exists a subgroup $H\sbs \G$ 
of finite index such that
the finite set $Y=\G/H$ carries a $\G$-invariant pseudometric extending the 
metric $\rho|K$, where $K$ is identified with its image in $Y$.
In virtue of lemma~\ref{l:appr}, the distance preserving map
$g\mapsto f(g)p$ from $(K,\mu)$ to $\U$ can be approximated by 
a distance preserving map $j:(K,\rho)\to \U$ such that $j(1)=f(1)p=p$
and 
\[
d(j(g), f(g)p)\le \diam (K,\e\lambda)=\e
\]
for every $g\in K$. 
Using the injectivity of $\U$, extend $j$ to a distance-preserving map 
(which we still denote by $j$) $Y\to \U$.
The natural action by isometries of $\G$ on $Y$ gives rise to an action by isometries
of $\G$ on $j(Y)$ such that $j$ becomes a morphism of $\G$-spaces. In virtue of 
Lemma~\ref{l:g-embed}, there exists a homomorphism $k:\Iso(j(Y))\to P=\Iso(\U)$ which
sends each isometry of $j(Y)$ to its extension over $\U$. Let $f': \G\to P$
be the composition of the homomorphism $\G\to \Iso(j(Y))$ corresponding to the
action of $\G$ on $j(Y)$, with $k$. Since the group $\Iso(j(Y))$
is finite, $f'$ has a finite range. We claim that $f'\in O$. It suffices to verify
that $f'(g)p$ and $f(g)p$ are $\e$-close with respect to $d$ for every $g\in K$. 
The $g$-shift $j(x)\mapsto j(gx)$ on
$j(Y)$ sends the point $p=j(1)$ to $j(g)$ which is $\e$-close to $f(g)p$. 
Since $f'(g)$ extends the $g$-shift on $j(Y)$, we have $f'(g)p=j(g)$, 
and thus $d(f'(g)p, f(g)p)=d(j(g),f(g)p)\le\e$. 

{\em Case 2.} If $A=\{p_1,\dots,p_n\}\sbs \U$
has more than one point, we replace $\G$ by the free
product $\G*F$, where $F$ is the free group on $n$ generators $b_1,\dots,b_n$, and use
a reduction to Case 1. Pick a point $p\in \U$, and let $F$ act
by isometries on $\U$ in such a way that the generator $b_i$ sends $p$ to $p_i$, 
$1\le i\le n$. Use this action to extend in an obvious way the homomorphism
$f:\G\to P$ to a homomorphism (still denoted by $f$) $\G*F\to P$. Free groups are
residually finite, and the free product of residually finite groups is residually
finite (cf. \cite{Gru}), so we can apply the result of the preceding
paragraph to the group $\G*F$. Consider the finite subset 
$L=\{gb_i: g\in K, 1\le i\le n\}$ of $\G*F$. According to Case 1,
there exists a homomorphism $f':\G*F\to P$ such that $f'$ has a finite range and
$f'(g)p$ is $\e/2$-close to $f(g)p$ for every $g\in L$. In particular, 
$t_i=f'(b_i)p$ is $\e/2$-close to $f(b_i)p=p_i$.
For every $g\in K$ and every $i=1, \dots, n$ 
the point $f'(g)p_i$ is $\e/2$-close to $f'(g)t_i$, while 
$f'(g)t_i=f'(g)f'(b_i)p=f'(gb_i)p$ is $\e/2$-close to 
$f(gb_i)p=f(g)f(b_i)p=f(g)p_i$. It follows that $f'(g)p_i$ and $f(g)p_i$
are $\e$-close.
Thus the restriction of $f'$ to $\G$ has the required properties. \qed

\subsection*{Proof of Corollary \ref{c:kirch}}
The group $\G=F_\infty\times F_\infty$ is residually finite. According to 
Theorem~\ref{th:app}, representations with a finite range form a dense
subspace of $\Rep(\G, \U)$. This is stronger than the Kirchberg property
$\circledast_3$. 
\qed

\section{Connes' Embedding conjecture and Kircherg's property of topological groups\label{s:cec}}

In this section we discuss 
Connes' Embedding Conjecture (CEC), which is
presently one of the main open problems in the theory of operator algebras.
As mentioned in Section~1, CEC is equivalent to the
conjecture that the unitary group $U_s(\ell^2)$ satisfies the Kirchberg property.
Thus our Corollary~\ref{c:kirch} can be viewed 
as an analogue of CEC for the Urysohn space and its group of isometries. 

Since this paper is aimed at experts in topological
groups as much as at those in operator algebras, for the benefit of the former we begin by recalling some basic facts of the theory of operator algebras,
referring the reader for a more detailed treatment, for instance, to
the books by Sakai \cite{Sa} and Takesaki \cite{Tak}. 

Recall that a {\em von Neumann algebra} $M$ is a unital $C^\ast$-algebra
which, regarded as a Banach space, is a dual space: there is a
(necessarily unique) Banach space $M_\ast$, the \textit{predual} of $M$,
with the property that $M$ is isometrically isomorphic to $(M_\ast)^\ast$. 
A von Neumann algebra with a separable predual is called {\em hyperfinite}
if it is generated, as a von Neumann algebra, by an increasing sequence of
finite-dimensional subalgebras. 

A von Neumann algebra $M$ is 
called a {\em factor} if the centre of $M$ is
trivial, that is, consists of scalar multiples of $1$. 
For example, the von Neumann algebra
${\mathcal L}(\ell^2)$ of all bounded linear operators on the Hilbert space
$\ell^2$ is a hyperfinite factor. 

Let $E_\alpha,\alpha\in A$ be a family of normed spaces, and let
$\xi$ be an ultrafilter on the index set $A$. The ({\em Banach space})
{\em ultraproduct} of
the family $(E_\alpha)$ along the ultrafilter $\xi$ is the linear space
quotient of the $\ell^\infty$-type direct sum
$E=\oplus^{\ell^\infty}_{\alpha\in A}E_\alpha$ by the ideal ${\mathcal I}_\xi$
formed by all collections $(x_\alpha)_{\alpha\in A}\in E$ with the property
\[\lim_{\alpha\to\xi}x_{\alpha} =0,\]
equipped with the norm
\[\norm x = \lim_{\alpha\to\xi}x_{\alpha},\]
where $(x_\alpha)$ is any representative of the equivalence class $x$.
If the ultrafilter $\xi$ is free, the ultraproduct is always a Banach space.
For a general theory of ultraproducts of normed spaces 
(also known in nonstandard analysis as \textit{nonstandard hulls}), see 
\cite{HI}.

The ultrapower of a family of $C^\ast$-algebras is again a $C^\ast$-algebra, but the property of being a factor is not necessarily preserved.
However, in the particular case where
all factors in a family are the so-called finite factors,
one can modify the construction of an ultraproduct so as to obtain a factor.

Recall that a (finite) \textit{trace} on a von Neumann algebra $M$ is a positive linear
functional $\tau\colon M\to\C$ with the property $\tau(AB)=\tau(BA)$ for
all $A,B\in M$. A trace $\tau$ is \textit{normalized} if $\tau(1)=1$.
One says that a factor $M$ is \textit{finite} if it admits a trace. One can show
that in this case the normalized trace on $M$ is unique. 
Finite factors of finite dimension are exactly all matrix algebras of the
form $M_n(\C)$, $n\in\N$. However, there exist finite factors that are
infinite-dimensional as normed spaces. They are called
\textit{factors of type $II_1$}. 

An example is
given by the following construction. Let $G$ be a (countable) discrete group.
Denote by $VN(G)$ the strongly closed unital $\ast$-subalgebra of 
${\mathcal L}(\ell^2(G))$ generated by all operators of left translation by
elements of $G$. This is the so-called ({\em reduced})
\textit{group von Neumann algebra}
of $G$. If all conjugacy classes of $G$ except for that of unity are
infinite, then $VN(G)$ is a factor of type $II_1$. For example, this is the 
case where $G=F_2$, the free group on two generators. On the contrary,
the factor ${\mathcal L}(\ell^2)$ does not admit a trace.

As was shown by Murray and von Neumann, there exists only one, up to an
isomorphism, hyperfinite factor of type $II_1$, denoted by $R$. For instance,
$R$ is isomorphic to the group von Neumann algebra of a locally finite
group (the union of an increasing sequence of finite subgroups) with
infinite conjugacy classes. (Probably the simplest example is $S^{fin}_\infty$,
the group of all permutations of the set $\N$ of natural numbers having a
finite support each.) 

Now let $M_\alpha$ be a family of finite factors, each equipped
with a normalized trace
$\tau_\alpha$, and let $\xi$ be an ultrafilter on the index set $A$.
The formula
\[\tau((x_\alpha)) = \lim_{\alpha\to\xi}\tau_\alpha(x_\alpha)\]
determines a trace on the Banach space ultraproduct $M$
of the family $(M_\alpha)$
along $\xi$. The subset 
\[{\mathscr I}_\xi = \{x\in M\colon \tau(x^\ast x)=0\}\]
is an ideal of $M$, and the factor-algebra $M/{\mathscr I}_{\xi}$ happens to be a finite factor,
called the \textit{von Neumann ultraproduct} of the family $(M_\alpha)$.
Under an obvious non-degeneracy assumption
(for every $n\in\N$, the set $\{\alpha\in A\colon \dim(M_\alpha)\geq n\}$
is in $\xi$),
the von Neumann ultraproduct $M/{\mathscr I}_{\xi}$ is non-separable, thus
has infinite dimension and is 
a factor of type $II_1$. For instance, the von Neumann
ultraproduct of all matrix algebras
$M_n(\C)$, $n\in\N$, equipped with their standard normalized traces, along any
free ultrafilter on the natural numbers, is a factor of
type $II_1$.

As every subfactor of a factor of type $II_1$ is again of type $II_1$,
one may wonder how large is the class of all separable factors of type $II_1$ embeddable
into ultrapowers of $R$. Such factors do not need to be hyperfinite:
already Connes had remarked \cite{connes-injective} that
$VN(F_2)$ is among them.

The following conjecture was formulated by Connes
in the same paper \cite{connes-injective}
(p. 105, third paragraph from the bottom).
\vskip 2mm

\noindent
{\bf Connes' Embedding Conjecture.}
{\it
Every factor of type $II_1$ embeds into an ultrapower of the
hyperfinite factor $R$ of type $II_1$. 
}\vskip 2mm

In the above conjecture, one can assume without loss in generality
that the factors have separable preduals, and the index set supporting the
ultrafilter is countable. Furthermore, one can replace
the ultrapower of $R$ with the von Neumann
ultraproduct of matrix algebras $M_n(\C)$. For a discussion, see e.g.
section 9.10 in \cite{pisier}.

In the last three decades, the conjecture has increased in importance 
and has become one of the main open problems of operator algebras theory.
Largely through the work of E. Kirchberg, numerous equivalent forms 
of Connes' conjecture came into existence.

By a {\em $C^\ast$-norm} on $*$-algebra we mean
a norm {\it the completion} with respect to which is a $C^\ast$-algebra.
If $A$ and $B$ are two unital $C^\ast$-algebras,
their algebraic tensor product $A\otimes B$ is not, in general, a
$C^\ast$-algebra again (unless one of the algebras is finite-dimensional), 
but it always 
admits at least one $C^\ast$-norm extending the norms on $A$ and $B$
(we identify $A$ and $B$ with their images under 
the natural embeddings $a\mapsto a\otimes 1$, $b\mapsto 1\otimes b$).
For instance, if $A$ is a $C^\ast$-subalgebra of
${\mathcal L}(\H_1)$ and $B$ is a a $C^\ast$-subalgebra of
${\mathcal L}(\H_2)$, then $A\otimes B$ embeds naturally as 
a $C^\ast$-subalgebra into
${\mathcal L}(\H_1\otimes\H_2)$ (the tensor product of Hilbert spaces),
and the norm induced by this embedding is called the {\em minimal tensor
product norm}. It has the remarkable property of being smaller than any
other $C^\ast$-norm on $A\otimes B$ (Takesaki). On the other hand, there exists
the {\em maximal tensor product norm,} which is the largest among all
$C^\ast$-norms on $A\otimes B$. The minimal and maximal tensor product norms
on $A\otimes B$ coincide in a number of important cases, for instance, when
one of the algebras $A,B$ is {\em nuclear}.

If no confusion can arise, one normally
denotes the $C^\ast$-algebra completion of
the algebraic tensor product $A\otimes B$, equipped with one or another norm,
with the same symbol $\otimes$, without adding a ``hat'' to mark the completion,
for instance, $A\otimes_{min} B$ is the \textit{completion} of the algebraic
tensor product equipped with the minimal tensor norm, and so forth.

If $A$ is a unital $C^\ast$-algebra, the {\em unitary group} of $A$ is a
multiplicative subgroup consisting of all unitaries of $A$, that is,
$u\in A$ with $u^\ast u=uu^\ast = 1$.
Every discrete group $G$ admits a universal embedding, as a subgroup,
into the unitary group of a suitable $C^\ast$-algebra. Namely, there exist a
unital $C^\ast$-algebra $C^\ast(G)$, called the ({\em full}) 
{\em group $C^\ast$-algebra} of $G$, and a group homomorphism
(in fact, a monomorphism), $i$, from $G$ to the unitary
group $U(C^\ast(G))$ with the
property that, whenever $A$ is a unital $C^\ast$-algebra and
$f\colon G\to U(A)$ is a group homomorphism, there is a unique morphism of
$C^\ast$-algebras $\bar f\colon C^\ast(G)\to A$ with $\bar f\circ i=f$.
The $C^\ast$-algebra $C^\ast(G)$ is the $C^\ast$-envelope of the $*$-algebra
$l^1(G)$, it is unique up to an isomorphism for every discrete group $G$.

Here is a useful example to consider: the full $C^\ast$-algebra of the
direct product $G\times H$ of two groups is naturally isomorphic to the
{\em maximal} tensor product $C^\ast(G)\otimes_{max} C^\ast(H)$.

Below is a statement which is equivalent to the Connes' Embedding
Conjecture \cite{kirchberg}, see also \cite{ozawa} or \cite{pisier}, ch. 16.
\vskip 2mm

\noindent
{\bf Conjecture (Kirchberg)}
{\it
The tensor product of the 
group $C^\ast$-algebra $C^\ast(F_2)$ of the free group on two generators
with itself admits a unique $C^\ast$-algebra norm. (That is, the {\it max}
and {\it min} norms on $C^\ast(F_2)\otimes C^\ast(F_2)$ coincide.) 
}\vskip 2mm

A {\it representation} of a $C^\ast$-algebra $A$ in a Hilbert space $\H$ is
a $C^\ast$-algebra morphism $\pi\colon A\to{\mathcal{L}}(\H)$.
The {\it essential space} of a representation $\pi$ is 
the closure of $\pi(A)(\H)$ in $\H$. 
A representation $\pi$ is {\it degenerate} if its essential space is a proper
subspace of $\H$, and {\it finite-dimensional} if the essential space 
is finite-dimensional. A representation $\pi$ of a unital $C^\ast$-algebra
is {\it unital} if $\pi(1)=\I_{\H}$. 
If $A$ is a unital $C^\ast$-algebra, then a representation $\pi$ of $A$ is
unital if and only if it is non-degenerate. 

A $C^\ast$-algebra $A$ is called {\em residually finite-dimensional} 
({\it RFD}) if
it admits a separating family of finite-dimensional representations.
For instance, the full group $C^\ast$-algebra $C^\ast(F)$ of the non-abelian
free group (on any number of generators) is RFD,
this is a result by Choi \cite{choi}. Strictly speaking, the
finite-dimensionality of the algebra $\pi(A)$ is necessary, but not
sufficient, for $\pi$ to be finite-dimensional: the representation of
the one-dimensional $C^\ast$-algebra $\C$ in $\ell^2$ given by
$\pi(\lambda)=\lambda\I$ has all of $\ell^2$ as its essential space.

At the same time, a (unital) algebra $A$ is RFD if and
only if it admits a separating family of (unital) representations with 
finite-dimensional image, simply because every finite-dimensional algebra
admits a faithful finite-dimensional representation. 

It is not difficult to verify that
the {\em minimal} tensor product of two residually finite-dimensional
$C^\ast$-algebras is again residually finite-dimensional, and also that
if the {\em maximal} tensor product of two $C^\ast$-algebras is
residually finite-dimensional, then the maximal norm on the tensor product
coincides with the minimal one. These observations lead to the following
further reformulation of CEC, 
noted for example by Ozawa \cite{ozawa}, Prop. 3.19.
\vskip 2mm

\noindent
{\bf Conjecture} (equivalent to CEC).
{\it
The group $C^\ast$-algebra $C^\ast(F_2\times F_2)$ is
residually finite dimensional.
}\vskip 2mm

While the group $F_2\times F_2$ is of course residually finite, 
residual finiteness of a group $\G$ is in general insufficient for
the algebra $C^\ast(\G)$ to be RFD \cite{bekka}.

If $A$ is a $C^\ast$-algebra,
then $\Rep(A,\H)$ stands for the set of all
(degenerate and non-degenerate) representations of $A$ in $\H$. 
Following Exel and Loring \cite{el}, equip the
set $\Rep(A,\H)$ with the coarsest topology
making all the mappings of the form
\begin{equation}
\label{eq:el}
\Rep (A,\H)\ni \pi\mapsto \pi(x)(\xi)\in \H,~x\in A,~\xi\in \H
\end{equation}
continuous. Clearly, this topology is
inherited from
$C_p(A,{\mathcal{B}}_s(\H))$; here the subscript ``{\it p}'' as usual,
stands for the topology of pointwise convergence, while
${\mathcal{B}}_s(\H)$ is the space ${\mathcal{B}}(\H)$ endowed with
the strong operator topology, that is, the topology induced from
$C_p(\H,\H)$. The basic 
neighbourhoods of an element $\pi\in\Rep(A,\H)$ are of
the form
\begin{eqnarray*}
\nonumber
{\mathcal{O}}_\pi[x_1,x_2,&\dots&,x_n;\Xi;\e] \\
&=&
\{\eta\in\Rep(A,\H)\colon
\norm{\pi(x_i)(\xi)-\eta(x_i)(\xi)}<\e,~i=1,2,\dots, n,~~\xi\in\Xi\},
\label{eq:elt}
\end{eqnarray*}
where $x_i\in A$ and $\Xi$ is a finite system of vectors in $\H$.

\begin{theorem}[Exel and Loring \cite{el}]
A $C^\ast$-algebra $A$ is residually finite-dimensional
if and only if the set of finite-dimensional
representations is everywhere dense in $\Rep(A,\H)$
for all Hilbert spaces $\H$. 
\label{th:al}
\end{theorem}


Every group $C^\ast$-algebra, $C^\ast(\G)$,
admits a \textit{counit}, that is, a one-dimensional unital
representation $\eta$, in $\H$, which is determined by the condition 
$\eta(g)=\Id$ for all $g\in \G$. 
Let $\pi\in\Rep(C^\ast(\G),\H)$. Denote $\H_1=\pi(C^\ast(\G))(\H)$ and associate
to $\pi$ the representation $\tilde\pi = (\pi\vert \H_1)\oplus 
\eta_{\H\ominus \H_1}$. This is a unital representation
of $C^\ast(\G)$ in $\H$, and if $\pi$ is finite-dimensional, then 
$\tilde\pi$ has a finite-dimensional image.
Notice that $\pi$ itself can be written in the form 
$\tilde\pi = (\pi\vert \H_1)\oplus {\mathbf{O}}_{\H\ominus \H_1}$.

For a unital $C^\ast$-algebra $A$, denote by $\Rep_1(A,\H)$ the subspace of
$\Rep(A,\H)$ consisting of unital representations.

\begin{corollary}
A group $C^\ast$-algebra $A=C^\ast(\G)$ is residually finite-dimensional
if and only if the set of unital representations with finite-dimensional image 
is everywhere dense in $\Rep_1(A,\H)$ for all Hilbert spaces $\H$. 
\end{corollary} 

\begin{proof} $\Rightarrow$: if a representation $\pi\in\Rep_1(A,\H)$ is
approximated by a net of finite-dimensional representations $(\pi_\alpha)$,
then $\pi$ is clearly approximated by the net $(\widetilde{\pi_\alpha})$ of
unital representations with finite-dimensional images. 
\par
$\Leftarrow$: let $\pi\in\Rep(A,\H)$ be arbitrary. We want to approximate
$\pi$ with finite-dimensional representations. Without loss in generality,
we may assume that $\pi$ is non-degenerate and so unital. 
There is a net $(\pi_\alpha)$ of unital
representations with finite-dimensional images approximating $\tilde\pi$. 
Let $x_1,x_2,\ldots,x_n\in A$, let $\Xi\in\H$ be finite, and let $\e>0$.
Find an $\alpha$ with 
\[
\norm{\pi(x_i)(\xi)-\pi_{\alpha}(x_i)(\xi)}<\e,
~i=1,2,\dots, n,~~\xi\in\Xi.
\]
Denote by $\H_1$ the linear subspace of $\H$ spanned by elements
$\pi_\alpha(x)(\xi)$, $x\in A$, $\xi\in\Xi$. This $\H_1$ is finite-dimensional
and invariant under all operators in $\pi_\alpha(A)$. The restriction
$\dot\pi_\alpha$ of $\pi_\alpha$ to $\H_1$ is a finite-dimensional 
representation of $A$, and 
\[
\norm{\pi(x_i)(\xi)-\dot\pi_{\alpha}(x_i)(\xi)}<\e,
~i=1,2,\dots, n,~~\xi\in\Xi.
\]
\end{proof} 

Unital representations of the 
group $C^\ast$-algebra $C^\ast(G)$ in a Hilbert space $\H$
are in a natural one-to-one correspondence with the unitary representations
of the group $G$ in $\H$, that is, group homomorphisms from $G$ to the
unitary group $U(\H)$ of the Hilbert space $\H$. We will equip the latter
group with the
\textit{strong operator topology}. This is simply the topology
of pointwise convergence on $\H$ (or on the unit ball), that is, a
topology induced by the embedding $U(\H)\subseteq C_p(\H,\H)$.
This topology makes $U(\H)$ into a Polish topological group. 

In view of the above remarks, there is a canonical
bijection $\Rep_1(C^\ast(G),\H)\leftrightarrow \Rep(G,\H)$, which
is in fact a homeomorphism.

The image $\pi(G)$ of a representation $\pi$ of a group $G$ is a
topological subgroup of the unitary group $U(\H)$ with the strong operator
topology. This image is a relatively compact subgroup if and only if 
$\pi$ factors through 
a strongly continuous representation of a compact group. 
For instance, this is the case where 
$\pi$ is the restriction to $G$ of 
a representation of $C^\ast(G)$ with a finite-dimensional image. 
Since every strongly continuous representation of a compact group, on the other hand, decomposes into a direct sum of finite-dimensional representations, one can easily deduce the following.

\begin{corollary}
For a discrete group $G$, the following conditions are equivalent:
\begin{enumerate}
\item the full $C^\ast$-algebra $C^\ast(G)$ is residually finite-dimensional;
\item if $\H$ is a Hilbert space, representations with 
relatively compact image (in the strong topology)
are everywhere dense in the space
$\Rep(G,\H)$. 
\end{enumerate}
\label{c:el}
\end{corollary}

Here it is enough to take a Hilbert space $\H$ of the same density character
as the cardinality of $G$. 

The topology considered by Exel and Loring
is finer than the well-known \textit{Fell topology}
\cite{fell}. For the Fell topology, an analogue
of the above characterization can be found in \cite{dixmier}.

We now summarize our discussion of various equivalent forms of CEC:

\begin{theorem}
\label{th:sum}
Each of the following conjectures is equivalent to Connes' Embedding
Conjecture:
\begin{enumerate}
\item the algebra $C^\ast(F_\infty\times F_\infty)$ is residually finite-dimensional;
\item the algebra $C^\ast(F_2\times F_2)$ is residually finite-dimensional;
\item the unitary group $U_s(\ell^2)$ has the Kirchberg property.
\end{enumerate}
\end{theorem}

The equivalence of $(1)$ and $(2)$ between themselves follows from the fact
that the algebras $C^\ast(F_\infty\times F_\infty)$ and $C^\ast(F_2\times F_2)$
are unital $C^\ast$-subalgebras of each other. The equivalence of $(1)$ and $(2)$ to Connes' Embedding Conjecture is Kirchberg's difficult result \cite{kirchberg}, cf. also \cite{ozawa,pisier}, while
the equivalence of $(1)$ and $(3)$ follows from  Corollary~\ref{c:el}, 
applied to the group $G=F_\infty\times F_\infty$.

Applying the property $\circledast_2$ with $H_2=\{e\}$, one gets the following
necessary condition for Kirchberg's property of a topological group $G$: every finite subset of $G$ can be
simultaneously approximated by elements of a finite set contained in a
compact subgroup. 
In particular, this is the case if $G$ admits an increasing chain of compact
subgroups with everywhere dense union.
This property is observed quite often in concrete ``large'' topological groups
of importance. Clearly, the unitary group $U(\ell^2)$ is one
of them. 
As shown by A.S. Kechris (private communication), an even stronger result
holds: every finite collection of elements of $U(\ell^2)$ is simultaneously
approximated with elements of a finite subgroup.
The group $\Iso(U)$ admits an increasing chain of finite subgroups with 
an everywhere dense union \cite{Ver05b}, cf. also \cite{P06}. 

Other Polish groups approximable by increasing chains of compact subgroups
include the infinite symmetric group $S_\infty$ of all self-bijections of
a countably infinite set and the group $\Aut(\I,\lambda)$ of measure-preserving
transformations of the standard Borel space with a non-atomic probability 
measure. For more examples see e.g. \cite{GP2,P05}. 

Theorem \ref{th:app} allows one to conclude that the Polish group
$\Iso(\U)$ has the Kirchberg property (Corollary \ref{c:kirch})
and leads us to 
ask the following:
does every Polish group $G$ approximated by an increasing chain of
compact subgroups satisfy Kirchberg's property?

If true, this will imply Connes' Embedding Conjecture when applied to
$G=U(\ell^2)$. At present, we are unaware of any topological group counterexamples.

\end{document}